\documentclass[12pt]{article}
\usepackage[right=1.25in,left=1.25in,top=1.1in,bottom=1.1in]{geometry}
\usepackage[utf8]{inputenc}
\usepackage{newunicodechar}
\DeclareUnicodeCharacter{1D54F}{\ensuremath{\mathbb{X}}}
\usepackage[font=small]{caption}

\usepackage{graphicx}
\usepackage{url}
\usepackage{amsmath}
\usepackage{amsthm}
\usepackage{engord}
\usepackage{float}
\usepackage{pdflscape}
\usepackage{booktabs}
\usepackage{pgfplots}
\pgfplotsset{compat=1.14}
\pgfplotsset{every axis label/.append style={font=\tiny}}
\usepackage[labelsep=period]{caption} 

\usepackage{amssymb} 
\usepackage{multirow} 

\usepackage{setspace}
\newcommand\fnsep{\textsuperscript{,}}
\newcommand\Prob{\mathbb{P}}

\usepackage{subcaption}
\usepackage{float}
\usepackage{natbib}
\usepackage[T1]{fontenc}    
\usepackage{hyperref}
\hypersetup{
	breaklinks=true,  
	colorlinks   = true, 
	urlcolor     = blue, 
	linkcolor    = blue, 
	citecolor   = blue 
}

\usepackage{etoolbox}
\makeatletter
\patchcmd{\NAT@citex}
{\@citea\NAT@hyper@{%
		\NAT@nmfmt{\NAT@nm}%
		\hyper@natlinkbreak{\NAT@aysep\NAT@spacechar}{\@citeb\@extra@b@citeb}%
		\NAT@date}}
{\@citea\NAT@nmfmt{\NAT@nm}%
	\NAT@aysep\NAT@spacechar\NAT@hyper@{\NAT@date}}{}{}

\patchcmd{\NAT@citex}
{\@citea\NAT@hyper@{%
		\NAT@nmfmt{\NAT@nm}%
		\hyper@natlinkbreak{\NAT@spacechar\NAT@@open\if*#1*\else#1\NAT@spacechar\fi}%
		{\@citeb\@extra@b@citeb}%
		\NAT@date}}
{\@citea\NAT@nmfmt{\NAT@nm}%
	\NAT@spacechar\NAT@@open\if*#1*\else#1\NAT@spacechar\fi\NAT@hyper@{\NAT@date}}
{}{}
\makeatother

\usepackage{amsfonts}       
\usepackage{lipsum}

\usepackage{mathpazo}

\usepackage{caption}

\usepackage{adjustbox} 
\usepackage{xcolor} 
\title{ \vspace*{-2.5cm} \hspace*{-0.5cm}
Penney's Game Odds From No-Arbitrage}

\author{Joshua B. Miller\thanks{University of Melbourne, Department of Economics.
\href{mailto:joshua.benjamin.miller@gmail.com}{joshua.benjamin.miller@gmail.com}\\ ORCID: 0000-0001-8412-3689}\ \thanks{I am grateful to the editor and anonymous referees, as well as to Ivan Balbuzanov, Faisal Sohail, and Tom Wilkening, for their helpful comments and suggestions. Any errors or omissions are my own.}
}

\date{ \vspace*{0.5cm} November 05, 2025\footnote{
Previous Version: July 21, 2025 (First Version: April 23, 2019).
}\\
} 

\begin{document}

\bgroup
\let\footnoterule\relax

\begin{singlespace}
\maketitle
\begin{center}
	{\large \textbf{Accepted for publication in \emph{Theory and Decision}.}}
\end{center}

\begin{abstract}
    \noindent 
   
    Penney's game is a two-player zero-sum game in which each player chooses a three-flip pattern of heads and tails, and the winner is the player whose pattern occurs first in repeated tosses of a fair coin. Because the players choose sequentially, the second mover has the advantage. In fact, for any three-flip pattern, there is another three-flip pattern that is strictly more likely to occur first. This paper provides a novel no-arbitrage argument that generates the winning odds corresponding to any pair of distinct patterns.  The resulting formula is equivalent to that generated by Conway's ``leading number'' algorithm. The accompanying betting-odds intuition adds insight into why Conway's algorithm works. The proof is simple and easy to generalize to games involving more than two outcomes, unequal probabilities, and competing patterns of various lengths. Additional results on the expected duration of Penney's game are presented. 
\end{abstract}



\end{singlespace}
\thispagestyle{empty}

\clearpage
\egroup
\setcounter{page}{1}

\section{Introduction\label{sec:introduction}}

\noindent Penney’s game is a two-player zero-sum game in which each player selects one of the eight possible patterns of three consecutive coin flips---$HHH$, $HHT$, $HTH$, $HTT$, $THH$, $THT$, $TTH$, or $TTT$---with the winner being the player whose pattern occurs first in a sequence of fair coin tosses. For instance, if the players choose $HTH$ and $TTH$, respectively, and the first five coin flips yield $THTTH$, then $TTH$ is the winning pattern. Remarkably, when players choose sequentially, the second mover always has a strategic advantage: for any pattern selected by the first mover, the second mover can respond with a pattern that is strictly more likely to appear first in the sequence \citep{Penney--JRM--1969}.

The implied intransitivity---the fact that the three-flip patterns cannot be (weakly) ordered from best to worst, even though some patterns are better than others in a head-to-head race---is widely viewed as paradoxical, as if Nature itself were being gamed. Martin Gardner, in his widely-read {\it Scientific American} column ``Mathematical Games,'' wrote: ``\dots most mathematicians simply cannot believe it when they first hear of it\dots\ It is certainly one of the finest of all sucker bets'' \citep{Gardner--ScientificAmerican--1974}. 

A formula for the odds of winning in a two-pattern match-up remained elusive until mathematician John Conway devised what he called the ``leading number'' algorithm, though he did not provide a proof.  Gardner wrote of the algorithm: ``I have no idea why it works. It just cranks out the answer as if by magic, like so many of Conway's other algorithms'' \citep[p.~306]{Gardner2001}.

This paper presents a straightforward and intuitive proof of Conway’s algorithm, based on the economic principle of no-arbitrage.  Unlike earlier proofs, which rely on advanced methods, this approach provides a transparent intuition for why the resulting odds take the form they do.\footnote{The proofs of Conway's algorithm have used combinatorial methods \citep{GuibasOdlyzko--JCTA--1981}, martingale optional stopping \citep{Li--AnProb--1980}, Markov chain imbedding \citep*{GerberLi--SPA--1981}, generating functions \citep[pp.~401--410]{GrahamKnuthPatashnik1994}, as well as other approaches \citep{StefanovPakes--AnAProb--1997,pevzner_dna_1993}. \citet{Li--AnProb--1980} and \citet{GuibasOdlyzko--JCTA--1981} attribute the first proof to S.\ Collings and a generalization to J.G.\ Wendel, both unpublished manuscripts.}\fnsep\footnote{See \cite{Nickerson--UMAP--2007} for extensive coverage, including a simple intuition for how the second mover can construct an advantageous response, as well as additional history. For proofs relating to the second mover's best response, see \citet{felix_optimal_2006}, which builds on \citet{GuibasOdlyzko--JCTA--1981}.} The argument generalizes naturally to i.i.d.\ sequences with more than two outcomes, asymmetric probabilities, and games with multiple competing patterns of varying lengths.

To build intuition for this result, we begin with the simpler case of two-flip patterns, where direct enumeration first reveals the game’s counterintuitive structure. Suppose the first mover selects the pattern $HH$ and the second mover responds with $TH$. In this case, the moment a tail appears within the first two flips, the second mover is guaranteed to win eventually, with the complete sequence taking one of the forms $TH$, $HT^kH$, or $TT^kH$, where $T^k$ denotes $k \geq 1$ consecutive tails. The probability of such a win is therefore $3/4$---equivalently, three chances in favor for every one against. This corresponds precisely to the odds generated by Conway’s algorithm.\footnote{In Conway’s formulation, the \emph{leading number} $XY$ measures the overlap between the ``trailing'' flips of $X$ and the ``leading'' flips of $Y$.  For two-flip patterns $X,Y \in \{H,T\}^2$, it is defined as $XY := 2\,\delta_{x_1 y_1}\delta_{x_2 y_2} + \delta_{x_2 y_1}$, where $\delta_{xy}:=1$ if $x=y$ and $0$ otherwise. Conway's algorithm calculates the odds ratio as $\tfrac{AA-AB}{BB-BA} = \tfrac{3-0}{2-1}=3$, which matches the direct enumeration above.}

While restricting Penney’s game to two-flip patterns removes the second mover’s {\it strict} advantage, the simplified setting is sufficient to fully illustrate the no-arbitrage argument---an argument that scales naturally to longer patterns, where direct enumeration becomes infeasible.

Consider two players, Ann and Bob, who agree to take opposite sides of any fair \textdollar 1 bet on a single coin flip, which pays \textdollar 2 to the winner and \textdollar 0 to the loser.   More generally, a fair bet on any pattern of consecutive coin flips can be replicated through a sequence of fair bets on individual flips.  For example, to bet on a specific two-flip pattern---say, $HH$---Ann can bet \textdollar 1 with Bob that the first flip will be $H$; if she wins, she uses the \textdollar 2 payout to bet that the second flip will also be $H$.  This {\it betting chain}, or ``$HH$-bet'', replicates a \textdollar 1 wager on the pattern $HH$ occurring in the next two flips, which, to preserve fairness, must pay \textdollar 4 if successful and \textdollar 0 otherwise. From Bob’s perspective, he is effectively selling Ann a derivative bet on the $HH$ pattern.

\begin{figure}[t]
  \centering
  \captionsetup{width=0.95\textwidth}
 \includegraphics[width=.95\textwidth]{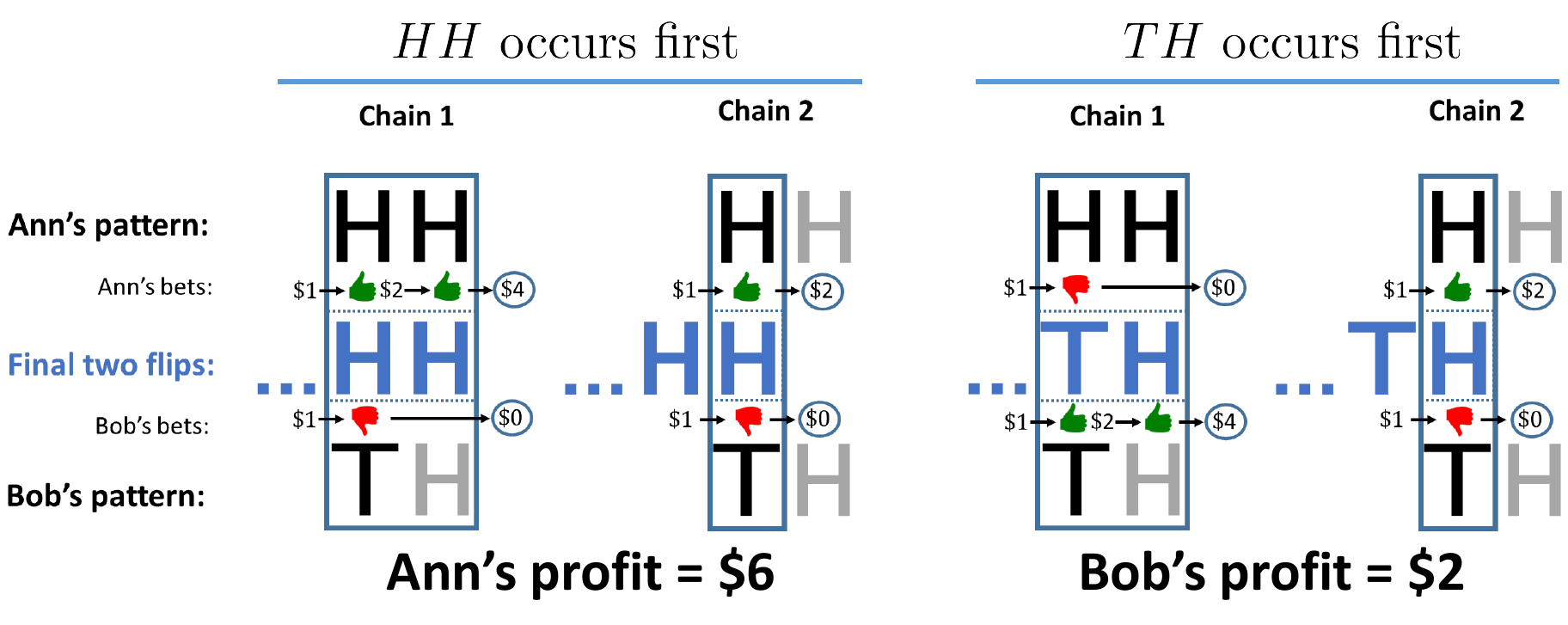}
  \caption{In the left panel, Ann's pattern $HH$ occurs before Bob's pattern $TH$; Ann's overall profit (from Bob) is \textdollar 6, gaining $\$ 4 +\$2 =\$6$ from Bob with her two $HH$ betting chains with him, and losing nothing to Bob from his two $TH$ betting chains with her.  In the right panel, pattern $TH$ occurs before pattern $HH$; Bob's overall profit (from Ann) is \textdollar 2, gaining $\$ 4 +\$0 =\$4$ from Ann with his two $TH$ betting chains with her, and losing $\$2$ to Ann from the first bet in her second $HH$ betting chain with him.  Overall, Bob has risked $\$6$ on $TH$ occurring first, and stands to receive a profit of $\$2$ if he is correct,  i.e., his betting odds are $\$3$ in favor to every $\$1$ against. If Bob is risk neutral, this implies a probability of $3/4$.
    } 
  \label{fig: Trading Strategy}
\end{figure}

In an ongoing sequence of coin flips, Bob can bet that his pattern---say, $TH$---occurs before Ann’s pattern, $HH$, by ``shorting'' her pattern to fund a ``long'' position on his own. To implement this strategy, before {\it each flip}, he sells a new $HH$-bet to Ann for \textdollar 1 and uses the proceeds to buy a $TH$-bet from her. Once one of the two patterns occurs, they agree to cease betting and settle the payoff of any incomplete betting chain.  If $TH$ occurs before $HH$, Bob earns an overall net gain of \textdollar 2, as illustrated in the right panel of Figure~\ref{fig: Trading Strategy}. From his long positions, Bob gains \textdollar 4 from the $TH$-bet he initiates on the penultimate flip and \textdollar 0 from the $TH$-bet he initiates on the final flip, due to the second chain's unsuccessful first-$T$ wager. From his short positions, he loses \textdollar 0 on Ann's $HH$-bet initiated on the penultimate flip and \textdollar 2 on the $HH$-bet she initiates on the final flip, due to the incomplete chain's successful first-$H$ wager. If instead $HH$ occurs before $TH$, Bob incurs an overall net loss of \textdollar 6 to Ann, as illustrated in the left panel of  Figure~\ref{fig: Trading Strategy}. He gains nothing from his long $TH$ positions, because the first-$T$ wager in each $TH$-bet never matches a flip, and he loses \textdollar 6 on his short $HH$ positions: \textdollar 4 to Ann's $HH$-bet initiated on the penultimate flip, and \textdollar 2 to her first-$H$ wager in the chain initiated on the final flip.

By using the proceeds from his short positions to finance his long positions, Bob has constructed a self-financing trading strategy in which he gains \textdollar 2 if $TH$ occurs before $HH$, and loses \textdollar 6 otherwise. As both Ann and Bob have accepted fair odds on the underlying bets, Bob cannot expect to profit from his strategy. Otherwise, he would have a successful statistical arbitrage opportunity---i.e., if  Bob's expected profits were positive, then he could renew his zero-cost self-financing trading strategy upon each occurrence of $TH$ or $HH$ in the ongoing sequence of flips. With positive expected profit each time he renews, he would attain any amount of wealth with probability one (in the limit).\footnote{Another way of demonstrating the necessity of this zero expected profit condition is to observe that Bob has constructed a ``gambling policy'' consisting of a ``betting system'' (long $TH$, short $HH$) and a stopping rule (stop at $TH$ or $HH$). For any gambling policy composed of fair bets, the expected fortune produced when betting stops must equal the initial fortune \citep[pp.99--100]{Billingsley1995}. In Bob's case, his initial fortune is \textdollar 0 \label{fn: Billingsley}.}

It follows directly from the zero expected profit condition for fair bets that the betting odds must equal the objective (i.e., probabilistic) odds:

\begin{align}
    &\underbrace{2\times\mathbb{P}(TH\prec HH) - 6\times\mathbb{P}(HH\prec TH)}_{\text{Bob's expected profit}} = 0 \nonumber \\
    &\hspace{6cm} \Longrightarrow \quad 
    \underbrace{\frac{\mathbb{P}(TH\prec HH)}{\mathbb{P}(HH\prec TH)}}_{\substack{\text{Objective} \\ \text{Odds}}} 
   \ \ = \
    \underbrace{\frac{6}{2} \vphantom{\dfrac{\mathbb{P}(TH\prec HH)}{\mathbb{P}(HH\prec TH)}}}_{\substack{\text{Betting} \\ \text{Odds}}}
    \label{eq:zero_profit}
\end{align}
where $\mathbb{P}(TH\prec HH)$ is the probability that pattern $TH$ occurs before pattern $HH$. Because Bob has bet \textdollar 6 in favor of $TH\prec HH$ for every \textdollar 2 Ann bet against it and each views their amount at stake as fair, the probabilistic odds must match: six chances in favor of $TH\prec HH$ for every two against. With eight chances overall, this implies $\mathbb{P}(TH\prec HH) = 6/8 = 3/4$,  consistent with the result obtained by direct enumeration.

We have, in effect, derived the risk-neutral probability that $TH$ occurs before $HH$ from the fact that the associated bet is fairly priced \citep[see, e.g.,][pp.61--63]{LeRoyWerner2003}.\footnote{While the existence of positive expected profits on Bob's trading strategy would not create an opening for a classic (risk-free) arbitrage opportunity, we have nevertheless derived the \emph{risk-neutral} probabilities associated with Bob's strategy in the following sense: Bob's gross return $r$ satisfies $r=8/6$ if $TH\prec HH$ and $0$ otherwise, while the risk-free (gross) return is $\bar{r}=1$. Therefore, the risk-neutral probabilities are those that satisfy the condition $\mathbb{E}[r]=\bar{r}$---that is, $\mathbb{P}(TH\prec HH)\cdot\frac{8}{6}=1$ \citep[see, e.g,][pp.61--63]{LeRoyWerner2003}.} It is straightforward to demonstrate that the resulting betting odds formula is identical to Conway’s leading number algorithm when initial bets are set at $\$0.50$ instead of $\$1$, which is shown in Section~\ref{sec: ConwayPenney}.

The betting odds approach illustrates how a core economic concept can illuminate problems in seemingly unrelated areas. In particular, the argument leverages the fact that a {\it statistical} arbitrage opportunity---a sequence of transactions that, in the limit, produces an investment that is costless, riskless, and has strictly positive profits---cannot arise from any trading strategy involving fairly priced financial assets.\footnote{The term ``statistical arbitrage'' comes from \citet{HoganJarrowTeoWarachka--JFE--2004}. The use of the term here is slightly different and is equivalent to \citet{Ross--JET--1976}'s limiting arbitrage opportunity \citep[for a discussion, see ][p.1287]{ChamberlainRothschild--Econometrica--1983}. In this setting, this is equivalent to the fact that fair games preclude successful gambling policies (see footnote~\ref{fn: Billingsley}).}

It is natural to ask whether the same betting logic provides an intuition that extends to the intransitivity of Penney’s classic game. It does, though the argument is more subtle. The details are developed in Section~\ref{sec: intransitive}, but the key insight follows from Equation~\ref{eq:zero_profit} and the two-flip example above: to obtain a probabilistic advantage under fair betting odds, Bob must choose a pattern $B\in\{T,H\}^2$ such that he risks losing more to Ann than Ann risks losing to him. Because Ann's pattern is given, Bob's potential loss on her $HH$-bets is fixed at $\$6$ when her pattern occurs first. Bob can then maximally reduce Ann's potential loss when his pattern occurs first by choosing $B$ to maximize her gains on $HH$-bets in that contingency. This occurs when the trailing flip of $B$ matches the leading flip of $HH$, i.e.\ when $B=bH$ with $b\in\{T,H\}$, which yields Ann a $\$2$ gain from her final $HH$-bet if Bob's pattern occurs first. Bob can then choose the first flip of his pattern to minimize his own gain from $bH$-bets when Ann's pattern occurs first, which occurs when the leading flip of Bob's pattern does not overlap with the trailing flip of Ann's pattern, i.e.\ when $B=TH$. This choice cements Bob’s probabilistic advantage, as it also minimizes his gain on his $bH$-bets when his pattern occurs first by ensuring that his second $bH$-bet cannot pay off.

While this betting logic illustrates a second-mover advantage in Penney’s game, it relies on the first-mover accepting even-money odds on the \emph{patterns themselves}, rather than on the underlying coin flips. In a competitive market, such an advantage would vanish as betting odds adjust to reflect the \emph{fair} odds implied by the underlying flips. The no-arbitrage principle---common to both financial and betting markets---ensures that any apparent advantage arising from an increase in the \emph{probability} of winning is precisely offset by a corresponding decrease in the \emph{payoff}, keeping expected net returns at zero.  In this sense, Nature remains immune to being ``gamed.''

\subsubsection*{Related literature and the limits to no-arbitrage}

The approach taken here to derive betting odds can be viewed more generally as an application of the optional stopping theorem. Under suitable regularity conditions, if $\tau\geq 1$ is a {\it stopping time} for a stochastic process $\{\pi_t\}_{t=0}^\infty$, then the expected value of the process at the stopping time must equal its initial value, $E[\pi_\tau]=E[\pi_0]$, provided that $\{\pi_t\}$ is a martingale---i.e., $E[\pi_{t+1}| \pi_t]=\pi_t$.  The argument here implicitly applies the theorem in much the same way as the modern approach to the {\it gambler's ruin problem}, with the stopping time being the first time the gambler either reaches a target gain or loses the initial stake, each of which is known in advance. With known terminal payoffs, the payoffs' corresponding probabilities are determined by the theorem's initial value condition.\footnote{In particular, let $x_t\in\{T,H\}$ denote coin flips and let $A,B\in\{T,H\}^3$ be Ann's and Bob's respective patterns.  Define the stopping time $\tau:=\min\{t:x_{t-2}x_{t-1}x_t\in\{A,B\}\}$, which depends only on information available up to time $t$. Equivalently, $\tau$ can be defined as the time at which Ann's profit, $\pi^A_t$, first equals either her winning payoff or the negative of Bob's, i.e.,  $\tau=\min\{t: \pi_t^A\in\{\pi_{AB},-\pi_{BA}\}\}$. Ann's profit at time $t$ satisfies $\pi_t^A=R_{x_{t-2}x_{t-1}x_t}(A)-R_{x_{t-2}x_{t-1}x_t}(B)=-\pi_t^B$. It is straightforward to show that the respective components are martingales, which implies that $\{\pi_t^A\}$ is a martingale. The probabilities for the winning payoffs follow directly from the optional stopping theorem's initial value condition, $E[\pi_\tau^A]=\pi_0^A=0$. See \citet[pp.~491-492]{grimmett_probability_2001} for a modern application of the optional stopping theorem to the gambler's ruin problem in the context of a random walk. For classical approaches see \citet[pp.~39--48]{Gorroochurn--book--2012}.\label{fn: optionalstop and Gf}}  The betting chain used in the argument is adapted from \citet{Li--AnProb--1980}’s proof of Solov’ev’s~\citeyearpar{Solovev--TPA--1966} formula, which shows, for example, why the expected waiting time before observing the pattern $HT$ is four flips, compared with six for $HH$---values that, in the betting framework above, equal the respective gains from going long on each pattern once it occurs.\footnote{In particular,  Li's~\citeyearpar{Li--AnProb--1980} proof uses a martingale argument to show that if $\tau_A$ is a stopping time for the first occurrence of pattern, say, $A\in\{T,H\}^3$, then $E[\tau_A]=R_A(A)$.  This can be briefly demonstrated in this setting by supposing that Ann goes long on $A$ each period, without any short positions. Following the approach from footnote~\ref{fn: optionalstop and Gf}, Ann's profit at time $t$ is given by $\pi_t'=R_{x_{t-2}x_{t-1}x_t}(A)-t$, so her expected profit at the stopping time must equal her initial profit, $E[\pi_\tau']=R_A(A)-E[\tau_A]=0$, yielding the result. Li's~\citeyearpar{Li--AnProb--1980} approach has been extended to other waiting time problems \citep{pozdnyakov_waiting_2006}.}



As a coin-flip paradox, Penney’s game belongs to a broader class of puzzles associated with the so-called {\it overlapping words paradox}, a term from DNA statistics, which refers to the fact that while the expected number of occurrences of a pattern in a sequence of random letters (each equally likely) depends only on its length, the distribution can differ significantly depending on how the pattern overlaps with itself \citep{pevzner_dna_1993}.  A striking example is {\it Litt's Game}, a recently developed two-player game in which each player selects a two-flip pattern $A,B\in\{T,H\}^2$, and after 100 flips the winner is whoever's pattern occurs most often \citep{litt_httpsxcomlittmathstatus1769044719034647001_2024}.   In contrast to Penney’s game---where patterns differing only in their final flip are evenly matched---in Litt’s game they are not. For instance, many are surprised to learn that $HT$ has a probabilistic edge over $HH$, despite both patterns having the same expected frequency \citep{segert_proof_2024, grimmett_alice_2025,janson_generalized_2025}. Additionally, Litt’s game is transitive, even when generalized to longer competing patterns, because each pattern’s position in the probabilistic ranking is determined solely by its degree of self-overlap---for example, among three-flip patterns, none ranks strictly higher than $HHT$, which has the lowest possible self-overlap \citep{basdevant_cases_2025}. In fact, the surprising advantage of $HT$ over $HH$ in Litt’s game arises from the same selection bias that underlies a well-known statistical error pervasive in the ``hot hand fallacy'' literature \citep{MillerSanjurjo--Econometrica--2018,miller_bridge_2019,miller_is_2021,miller_cold_2024}.  In particular, within a finite sequence of coin flips, while there is a 50-50 chance of heads on any flip that happens to immediately follow a streak of heads (or not), this does not describe the behavior of flips that are selected {\it because} they immediately follow a streak of heads.   Specifically, let $\hat{P}(H \mid H^k)$ denote the proportion of heads among those flips that are immediately preceded by at least $k$ consecutive heads, where $k \geq 1$. While it remains true that $\mathbb{P}(H|H^k)=1/2$ at any given flip within the sequence, the estimator for this conditional probability,  $\hat{P}(H|H^k)$, exhibits two counterintuitive properties: (i) its expected value is strictly less than $1/2$, and (ii)  its realized value is more likely to fall below $1/2$ than to exceed it.  For $k=1$, the latter property is equivalent to the advantage $HT$ has over $HH$ in Litt's Game.\footnote{To see this, let \#HH and \#HT represent the number of occurrences of the respective patterns. If at least one is positive, i.e., $\#H*:=\#HH+\#HT>0$, then $\#HH<\#HT\iff \hat{P}(H|H)=\frac{\#HH}{\#HT+\#HH}<\frac{1}{2}$. The identity, \[\frac{\mathbb{P}(\#HH<\#HT)}{\mathbb{P}(\#HH>\#HT)}=\frac{\mathbb{P}(\hat{P}(H|H)<\frac{1}{2}|\#H*>0)}{\mathbb{P}(\hat{P}(H|H)>\frac{1}{2}|  \#H*>0)},\]  which demonstrates the equivalence between the two relative odds ratios---and, therefore, the two problems---can be derived by multiplying and dividing the right-hand side of the identity by $\mathbb{P}(\#H*>0)$ and applying the law of total probability after including in the sum the respective complementary (null) probabilities, $\mathbb{P}(\#HH>\#HT, \#H*=0)=\mathbb{P}(\#HH<\#HT, \#H*=0)=0$. }\fnsep\footnote{\citet{MillerSanjurjo--Econometrica--2018} prove the first property, but not the second. While the second property is visually apparent in the histogram representing $\mathbb{P}(H|H^k)$'s sampling distribution, the published version of \citet[Fig.~3]{MillerSanjurjo--Econometrica--2018} includes only the histogram for $\hat{P}(H|H^k)-\hat{P}(H|T^k)$, from which we observe that it is more likely that $\hat{P}(H|H)<\hat{P}(H|T)$ than the reverse. This inequality is equivalent to  $\#HH \times\#TT<\#HT\times \#TH$.}

While the overlapping words paradox provides a common thread between these puzzles, the betting outcomes in Litt’s game (and the hot hand) differ in a way that exposes a limit to the no-arbitrage approach.  For example, in Litt's game, while a player can finance $HT$-bets by selling $HH$-bets, the associated stopping time is fixed and independent of any payoffs that can be known in advance. Because payoffs at the stopping time are random variables in this setting, probabilities cannot be derived from the optional stopping theorem's initial value (no-arbitrage) condition.\footnote{While the payoffs in the event a pattern wins are not pinned down, if the expected payoffs in each winning contingency could be determined, then the respective probabilities would be recoverable from the no-arbitrage condition.}

The remainder of this paper proceeds as follows: Section~\ref{sec: intransitive} demonstrates the intransitivity of Penney's classic game using the betting odds approach. Section~\ref{sec: Penney} derives a generalized version of Conway's formula using the no-arbitrage argument in a setting involving patterns of potentially different lengths and outcomes drawn from an arbitrary i.i.d.\ categorical distribution, and then establishes its connection to Conway's leading number algorithm.  Section~\ref{sec: Additional Results} extends this approach to waiting times and illustrates properties of the canonical three-flip Penney's game.\footnote{All source code used to validate the results is available at \href{https://github.com/joshua-benjamin-miller/penneysgame}{https://github.com/joshua-benjamin-miller/penneysgame}.
A rendered version of the Python notebook session can be found at:
\href{https://github.com/joshua-benjamin-miller/penneysgame/blob/master/Validating-Penney.ipynb}{https://github.com/joshua-benjamin-miller/penneysgame/blob/master/Validating-Penney.ipynb}}

\section{Why Penney’s Game is Intransitive \label{sec: intransitive}}

While the main contribution of the betting framework lies in its intuitive proof of Conway’s algorithm, it also provides a clear demonstration of the intransitive structure of Penney’s classic game, without requiring the computation of odds for every possible pair of competing patterns.

Suppose Ann and Bob choose competing three-flip patterns while agreeing to take opposite sides of fair bets on the coin flips underlying their chosen patterns, as outlined above. If Ann selects her pattern first, then Bob can ensure that his pattern is strictly more likely to appear before Ann's by choosing one in which his resulting (fair) stake---the amount he risks losing---is greater than Ann's.

Specifically, let $A=a_1a_2a_3\in\{T,H\}^3$ denote Ann's pattern.  The task is to construct $B\neq A$ such that Ann’s winning payoff (Bob's stake) is greater than Bob's winning payoff (Ann's stake).  Ann’s and Bob’s payoffs can be expressed in terms of their respective betting chains, ``$A$-bets'' and ``$B$-bets.'' If $A$ occurs before $B$, written $A\prec B$, then Ann's winning payoff is
\[\pi_{AB}=R_A(A)-R_A(B),\]
where $R_A(A)\in\{8,10,12,14\}$ is Ann's gross return from her three $A$-bets (the earliest of which must succeed), and $R_A(B)\in\{0,2,4,6\}$ is Bob's gross return arising from his three $B$-bets (the earliest of which cannot succeed for Bob).   Analogously, if $B\prec A$, then Bob's winning payoff is 
\[\pi_{BA}=R_B(B)-R_B(A),\]
where $R_B(B)\in\{8,10,12,14\}$ is Bob's gross return from his three $B$-bets (the earliest of which must succeed), and $R_B(A)\in\{0,2,4,6\}$ is Ann's gross return arising from her three $A$-bets (the earliest of which cannot succeed for Ann). Observe that Ann's winning payoff is strictly greater than Bob’s if and only if 
\begin{equation}
    R_A(A)+R_B(A)>R_B(B)+R_A(B), \label{eq:sum gains}
\end{equation}
that is, if, summed across both contingencies, Ann's potential gains on her $A$-bets exceed Bob's potential gains on his $B$-bets---equivalently, if Bob’s aggregate stake across both contingencies, on his pattern, exceeds that of Ann.


Given Ann's choice of A, Bob can maximize the left-hand side of \eqref{eq:sum gains} by focusing on $R_B(A)\in\{0,2,4,6\}$, Ann's gain if $B\prec A$, because he cannot affect $R_A(A)$, Ann's gain if $A\prec B$.  By choosing a pattern that ends with the leading two flips of Ann's pattern---i.e., choosing $B:=ba_1a_2$, $b\in\{T,H\}$---Bob ensures that if $B\prec A$, then Ann's earliest---i.e, highest-return---$A$-bet that is successful is her second one, which yields Ann $R_B(A)\geq 4$.  For his $B$-bets on the right-hand side of \eqref{eq:sum gains},  Bob can cap his gain if $A\prec B$ to $R_A(B)\leq 2$ by setting $B:=\tilde{a_2}a_1a_2$, which ensures that his second $B$-bet---the earliest one that could potentially succeed---does not pay off. With this choice, Bob's gain if $B\prec A$ is minimized to $R_B(B)=8$ because only his $B$-bet pays off.\footnote{His second chain, $\tilde{a_2}a_1$, fails to match the final two flips,  $a_1a_2$, regardless of whether $a_1=a_2$ or not. If $a_1=a_2$,  the first bet (on $\tilde{a_2}$) fails, whereas if $a_1\neq a_2$, the second bet (on $a_1$) fails.} 

Together, $R_A(A)\geq 8$ and $R_A(B)\leq 2$ imply that Ann's winning payoff if $A\prec B$ satisfies $\pi_{AB}=R_A(A)-R_A(B)\geq 6$ while $R_B(B)=8$ and $R_B(A)\geq 4$ imply that Bob's winning payoff if $A\prec B$ is strictly lower, with $\pi_{BA}=R_B(B)-R_B(A)\leq 4$. Therefore Bob has a strictly higher chance of winning, with the betting odds in his favor statisfying $\frac{\pi_{AB}}{\pi_{BA}}\geq \frac{6}{4}=\frac{3}{2}$, i.e., $\mathbb{P}(A\prec B)\geq \frac{3}{5}=0.6$. This demonstrates that Bob has a second-mover advantage in choosing the pattern most likely to occur first, i.e., that Penney's game is intransitive.\footnote{For longer competing patterns of identical length $n>3$, the pattern that makes Bob strictly more likely to win still ends with the leading flips of Ann’s pattern, $B=ba_1a_2\dots a_{n-1}$, following the same reasoning provided here that maximizes $R_A(A)+R_B(A)$. Characterizing Bob's first flip $b$ is less straightforward and involves a trade-off between capping Bob's long-$B$ payoffs if he wins, $R_B(B)$, vs. if Ann wins, $R_A(B)$. See \citet{GuibasOdlyzko--JCTA--1981} and \citet{felix_optimal_2006} for alternative approaches.}

\section{The General Case \label{sec: Penney}}
Penney's game can be generalized to arbitrary pairs of patterns and to any i.i.d.\ sequence of a categorical random variable. Let $X_t$, $t=1,2,\dots$ be a sequence of i.i.d.\ draws from the categorical distribution $\mathbb{P}(X_t=j)=p_j \in [0,1]$ over a set of $m$ characters indexed by $j\in J=\{1,\dots,m\}$. Let $A,B\in \bigcup_{k=1}^\infty J^k$, with $J^k:=\times_{i=1}^k J$, be two patterns of respective lengths $k$ and $k'$ where individual characters of the pattern are represented by $A=a_1a_2\cdots a_k$ and $B=b_1b_2\cdots b_{k'}$, respectively.   Assume that neither pattern is nested within the other and that the sequence stops at the first occurrence of one of the patterns; let $\tau$ refer to this stopping time.

Following the same procedure described in Section~\ref{sec:introduction}, Ann and Bob are both willing to take either side of a fair bet on individual draws $X_t$, where a \textdollar 1  bet on the outcome $j$ yields $\$ \frac{1}{p_j}$ if $j$ occurs, and $\$0$ otherwise. They bet each other on whether pattern $A$ or $B$ will appear first, with betting odds defined by the underlying fair bets. In particular, Bob shorts Ann's pattern, $A$, and uses the proceeds to go long on his pattern, $B$. To implement this, at each draw  Bob sells to Ann a \textdollar 1 fair bet on outcome $a_1$, where Ann's potential winnings this draw, $\$\frac{1}{a_1}$, will be placed on outcome $a_2$ at the next draw, and so on until either $A$ occurs, or her betting chain is bankrupt.  At each draw, Bob uses the proceeds from Ann to buy from her a \textdollar 1 fair bet on outcome $b_1$, where Bob's potential winnings on this draw, $\$\frac{1}{b_1}$, will be placed on outcome $b_2$ at the next draw, and so on until either $B$ occurs, or his betting chain is bankrupt. Let $A\prec B$ correspond to the event that pattern $A$ occurs before pattern $B$; this means that the final $k$ draws of the sequence satisfy $(X_{\tau-k+1},\dots, X_{\tau})=(a_1,\dots,a_k)$ and the fair (gross) return for a successful bet at each of these respective draws is given by $\left(\frac{1}{\mathbb{P}(X_t=a_{1})},\dots,\frac{1}{\mathbb{P}(X_t=a_{n})}\right)=\left(\frac{1}{p_{a_1}},\dots,\frac{1}{p_{a_n}}\right)$.  If pattern $A$ appears first at time $\tau$, the fair payoff to Bob for going long on pattern $B$ is equal to the sum of the payoffs generated from each reinvested betting chain on pattern $B$ initiated from draw $t={\tau-k+1}$ to draw $t=\tau$ that is not bankrupt at time $\tau$, i.e., each betting chain initiated at draw $t$ where $B_1^{\tau-t+1}:=(b_1,\dots, b_{\tau-t+1})=(X_t,\dots X_\tau)=(a_{k-(\tau-t)},\dots,a_k)=:A_{k-(\tau-t)}^k$, which corresponds to when the leading characters of $B$, $B_1^{\tau-t+1}$, match the trailing characters of $A$, $A_{k-(\tau-t)}^k$. The fair payoff from the \textdollar 1 staked on the betting chain initiated at time $t$ that is not bankrupt at time $\tau$ is the product of the fair payoff from bets on the individual draws,  $\frac{1}{p_{a_{k-(\tau-t)}}}\times \frac{1}{p_{a_{k-(\tau-t-1)}}}\times\cdots\times \frac{1}{p_{a_{k}}}$.  Because Bob's overall payoff when pattern $A$ appears first only depends on the final $k$ draws, his fair payoff from going long on pattern $B$ satisfies
\begin{equation}
    R_A(B):=\sum_{s=1}^{k}1_{[A_s^k=B_1^{k-s+1}]} \prod_{i=0}^{k-s+1} \frac{1}{p_{a_{s+i}}}=\sum_{s=1}^{k}\prod_{i=0}^{k-s+1} \frac{1}{p_{a_{s+i}}}1_{[a_{s+i}=b_{1+i}]}
\end{equation}\label{eqn: general payoff}
In the same way, the fair payoff to Ann if pattern $A$ occurs first is given by $R_A(A)$. On the other hand, if pattern $B$ occurs before pattern $A$, the fair payoff to Ann is $R_B(A)$, whereas the fair payoff to Bob is $R_B(B)$. From Bob's perspective, he risks losing $R_A(A)-R_A(B)$ to Ann when pattern $A$ appears first for the chance of gaining $R_B(B)-R_B(A)$ in the event pattern $B$ occurs first.  The zero expected profit argument presented in Section~\ref{sec:introduction} applies directly in this setting, and therefore Bob's expected profit is zero if and only if the betting odds reflect the probabilistic odds.  The amount Bob risks losing by betting on $B$ relative to his potential payoff constitutes the fair betting odds in favor of $B$ (against $A$) occurring first, $R_A(A)-R_A(B):R_B(B)-R_B(A)$, which yields the probabilistic odds of $B$ occurring before $A$:
\begin{equation}
    \frac{\mathbb{P}(B\prec A)}{\mathbb{P}(A\prec B)}=\frac{R_A(A)-R_A(B)}{R_B(B)-R_B(A)}
\end{equation}\label{eqn: gen odds}
Therefore, the probability of $B$ occurring before $A$ is equal to the dollars (chances) risked in its favor divided by the total dollars at stake, i.e.:
\begin{equation}
\mathbb{P}(B\prec A)=\frac{R_A(A)-R_A(B)}{R_A(A)-R_A(B)+R_B(B)-R_B(A)}
\label{eqn: gen prob}
\end{equation}
This game and argument can be extended to more than two competing patterns, with Bob going long on one pattern by using the proceeds from shorting the remaining patterns.\footnote{More precisely, the generalization to  $m>2$ competing patterns $P_j$, $j=1,\dots,m$, involves Bob going long on pattern $P_1$  by using the proceeds from shorting the remaining patterns $P_2,\dots, P_{m}$, with Ann taking the other side.  Bob then repeats this process, each time with a new (independent) sequence,  for $j=2,\dots, m-1$,  creating a position in which he goes long on pattern $P_j$ using the proceeds from shorting the remaining patterns $P_{j'}$ for $j'\neq j$, with Ann taking the other side.  Each of the $j=1,\dots,m-1$ independent positions will produce the zero expect profit equation
\[\Prob(P_j \text{ appears first})\times G_{P_j} +\sum_{j'\neq j} \Prob(P_{j'} \text{ appears first}) \times (- L_{P_{j'}})=0\]
where Bob's gain with $P_j$ appears first, $G_{P_j}$, is the gain from the long position that is financed from the $m-1$ short positions, i.e. $ G_{P_j}=(m-1)R_{P_j}(P_j)-\sum_{j'\neq j}R_{P_j}(P_{j'})$, while the loses are each given by  $L_{P_{j'}}=(m-1)R_{P_{j'}}(P_j)-\sum_{j''\neq j'}R_{P_{j'}}(P_{j''})$, for $j'\neq j''$. We can solve these equations for  $\Prob(P_{j} \text{ appears first})$ for each $j$ with  using the fact that the probabilities sum to $1$.}

\subsection{Relationship with Conway's Leading Number Algorithm\label{sec: ConwayPenney}}
Conway's leading number algorithm, while fully generalizable, is defined for the canonical three-flip Penney's game and produces the odds in favor of pattern $B$ (against pattern $A$):
\begin{equation}
   AA-AB:BB-BA  
\end{equation}\label{eqn: conway odds}
where for any two patterns $A=a_1a_2a_3$, $B=b_1b_2b_3$, the leading number $AB$ is a measure of how the leading characters of $B$ overlap with the trailing characters of $A$, and is given by the formula:
\begin{equation}
  AB:=2^2\prod_{i=1}^3 [a_i=b_i] + 2\prod_{i=1}^2 [a_{i+1}=b_i]+ 2^0[a_3=b_1]
\end{equation}\label{eqn: conway leading number}
The leading number $AB$ is tightly related to the payoff for pattern $B$ betting chains when $A$ occurs first, $R_A(B)$. Following Equation~\ref{eqn: general payoff}, when $A$ occurs first, this payoff satisfies $R_A(B)=2^3\prod_{i=1}^3 [a_i=b_i] + 2^2\prod_{i=1}^2 [a_{i+1}=b_i]+ 2[a_3=b_1]=2AB$. This means that Conway's algorithm corresponds to the betting chain beginning with $\$0.5$ a bet, instead of a $\$1$ bet. The betting payoff, $R_A(B)$, and Conway's leading number, $AB$, measure the same pattern overlap properties---the degree to which the leading flips of pattern $B$ overlap with the trailing flips of pattern $A$.    Therefore, while the betting odds $R_A(A)-R_A(B):R_B(B)-R_B(A)$ differ from Conway's odds formula $AA-AB:BB-BA $, this a matter of scaling, which does not affect the objective probabilistic odds, or probabilities.

\section{Additional Results\label{sec: Additional Results}}
In this section, we derive the expected time waiting for the game to end and explore how the first mover's disadvantage depends on the probabilities.
\subsection{Expected Length of Penney's Game}
The generalized leading number $R_A(B)$ can be used to calculate the expected waiting time until pattern $A$ or $B$ occurs for any two patterns $A,B\in \bigcup_{k=1}^\infty J^k$. The argument is a relatively straightforward adaptation and extension of an argument presented in \citet{Li--AnProb--1980}.

Suppose a Carl offers fair bets on a sequence $\{X_t\}$ of i.i.d.\ characters according to the above categorical distribution. Carl terminates betting the moment pattern $A$ or $B$ occurs. Further, assume that Ann and Bob place bets with Carl according to patterns $A$ and $B$, respectively. Below we see why the expected waiting time until pattern $A$ or $B$ occurs at trial $\tau$, is given by
\begin{equation}
E[\tau]=
\begin{cases}
R_A(A) & \text{if } A\subseteq B\\
R_B(B) & \text{if } B\subseteq A \\
\frac{R_A(A)\times R_B(B)-R_A(B)\times R_B(A)}{R_A(A)+R_B(B)-[R_A(B)+R_B(A)]} & \text{otherwise }
\end{cases}
\end{equation}
where $A\subseteq B$ represents the case in which $A$ is a substring of $B$. The terms $R_A(A)$ and $R_A(B)$ represent Ann and Bob's respective payoffs if pattern $A$ occurs first.

To see why this equation holds, note that when the Carl stops offering bets at time $\tau$, Ann and Bob will have each initiated $\tau$ bet sequences and will have paid the casino $\mathdollar\tau$ each. Let $R_{A\cup B}$ represent the total amount that Carl pays out to Ann and Bob when betting ends. Because Carl is offering a fair game, as argued in the previous section (see, for example, Footnote~\ref{fn: Billingsley}), Carl's expected profits must be equal to zero. With the expected profits given by $E[2\tau-R_{A\cup B}]$, this implies that $E[\tau]=E[R_{A\cup B}]/2$, by linearity of expectations.

In the case that one pattern is a substring of the other,  then betting will stop when the substring occurs, in which case Ann and Bob's bets perfectly overlap and the Carl pays an identical amount to Ann and Bob, which implies that if $A\subseteq B$, then $E[\tau]=E[R_{A\cup B}]/2=E[2R_A(A)]/2=R_A(A)$, and the expected waiting time is equal to the (known) payoff to Ann in this contingency. In the case that patterns are distinct and neither pattern is a substring of the other, we can derive a formula for the expected waiting time until pattern $A$ or $B$ occurs:

\begin{align}
E[\tau] &=\frac{E[R_{A\cup B}]}{2}\notag \\
        &=\frac{E[R_{A\cup B}|A\prec B]\mathbb{P}(A\prec B)+E[R_{A\cup B}|B \prec A]\mathbb{P}(B\prec A)}{2}\notag\\
        &=\frac{[R_A(A)+R_A(B)]\mathbb{P}(A\prec B)+[R_B(A)+R_B(B)]\mathbb{P}(B\prec A)}{2}\\
        &=\frac{R_A(A)\times R_B(B)-R_A(B)\times R_B(A)}{R_A(A)+R_B(B)-[R_A(B)+R_B(A)]}
\end{align}
Expression~(8) follows because in the case of event $A\prec B$, $R_A(A)$ is the amount paid to Ann, and $R_A(B)$ is the amount paid to Bob, etc.  Expression~(9) follows by combining Expression~(8) with Equation~\ref{eqn: gen prob}.

\subsection{Optimal First Mover Choices}
The resulting formulae can be used to investigate various aspects of Penney's game. For example, we can consider whether the first mover is always at a disadvantage.
\begin{figure}[t]
    \centering
    \begin{subfigure}[b]{0.45\textwidth}
              \captionsetup{width=\textwidth}
         \includegraphics[width=\textwidth]{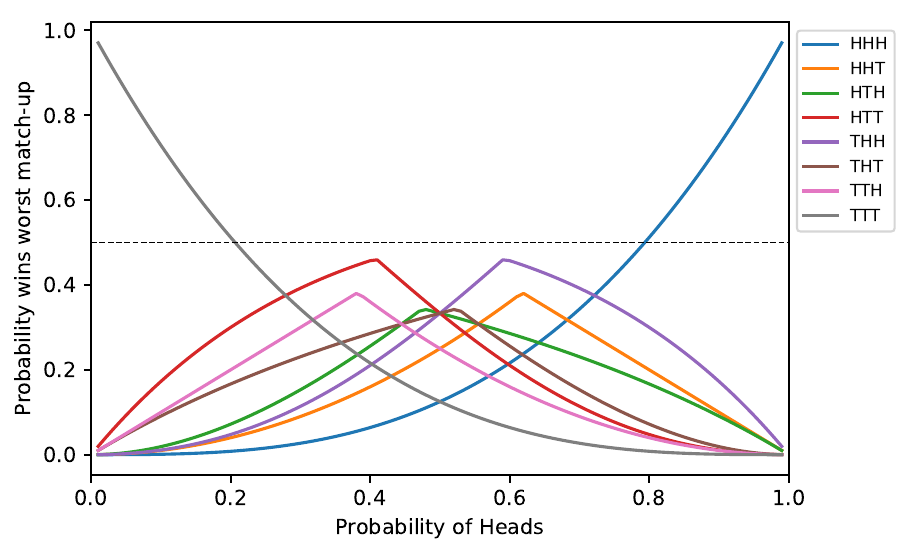}
        \caption{Probability pattern appears first}
        \label{fig: probability}
    \end{subfigure}
    ~ 
    \begin{subfigure}[b]{0.45\textwidth}
          \captionsetup{width=\textwidth}
         \includegraphics[width=\textwidth]{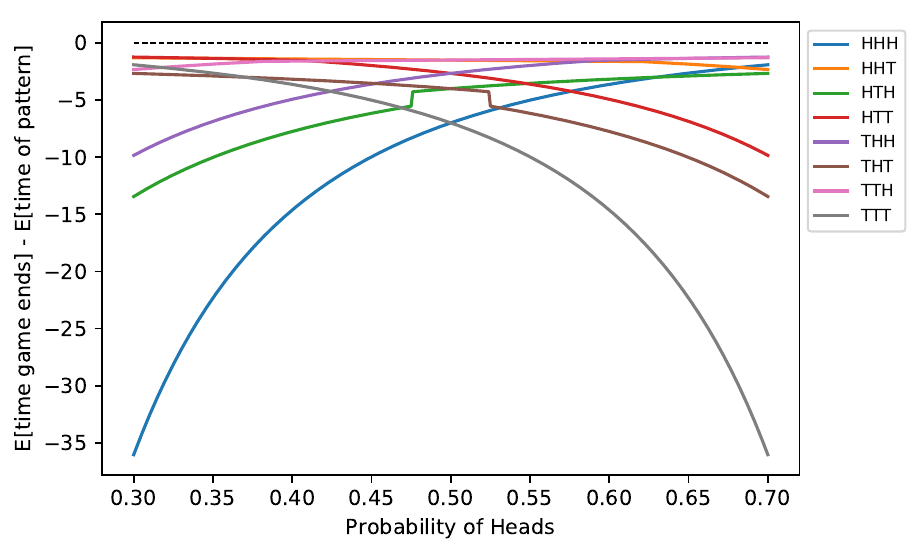}
        \caption{Change in waiting time}
        \label{fig: waiting time}
    \end{subfigure}
    \caption{Panel (a) the probability that the first mover's pattern appears before the second mover's best response is plotted as a function of the probability of heads; Panel (b) the decrease in expected time waiting for a pattern to occur when the second mover's pattern is added to the first mover's pattern, plotted as a function of the probability of heads.}\label{fig: worst matchup}
\end{figure}

Let Ann be the first mover and Bob be the second mover. Ann assumes that Bob will best respond to her choice, i.e.\, he will maximize the probability that his pattern appears first. In Figure~\ref{fig: probability}, the probability that Ann's pattern occurs before Bob's pattern is plotted as a function of the probability of heads.  As can be seen, when the probability of heads is sufficiently high or low, Ann gains the advantage in Penney's game.

Figure~\ref{fig: waiting time} reports the decrease in Ann's expected waiting time for the game to end vs. her expected waiting time for her chosen pattern. Unsurprisingly, the pattern that gives Ann the best chance of beating Bob's best response isn't likely to lead to a substantially shorter game when compared to the situation of Ann waiting alone for that pattern to occur.

\section{Conclusion\label{sec: Conclusion}}
In Penney's game, the objective odds of the event that pattern $A$, say $HTH$, occurs before pattern $B$, say $TTH$, has a simple representation using Conway's leading number algorithm. This suggests a simple explanation. We have shown how to construct a simple trading strategy involving fair bets on coin flips that amounts to a bet on the event that $A$ occurs before $B$. The absence of statistical arbitrage opportunities, particularly the impossibility of generating positive expected profits in a fair game, implies that the derived betting odds must equal the objective odds, yielding a simple formula for Penney's game odds. Notably, the betting-odds interpretation provided by the no-arbitrage argument offers novel insights into Penney's game and Conway's leading-number algorithm. Furthermore, the trading strategy can be readily generalized to games involving more than two outcomes, unequal probabilities, and competing patterns of various lengths. This proof illustrates how a core principle of economics---the idea that there is no such thing as a free lunch---can yield insight in unexpected areas.

\clearpage
\begin{singlespace}
\bibliographystyle{aer}
\bibliography{references.bib,zotero1.bib}
\end{singlespace}

\end{document}